\documentstyle[amscd,amssymb,verbatim,epsf]{amsart}

\begin{document}
\theoremstyle{plain}
\newtheorem{Thm}{Theorem}
\newtheorem{Cor}{Corollary}
\newtheorem{Ex}{Example}
\newtheorem{Con}{Conjecture}
\newtheorem{Main}{Main Theorem}
\newtheorem{Lem}{Lemma}
\newtheorem{Prop}{Proposition}

\theoremstyle{definition}
\newtheorem{Def}{Definition}
\newtheorem{Note}{Note}

\theoremstyle{remark}
\newtheorem{notation}{Notation}
\renewcommand{\thenotation}{}

\errorcontextlines=0
\numberwithin{equation}{section}
\renewcommand{\rm}{\normalshape}%

\title[Totally null surfaces]%
   {Totally Null Surfaces in Neutral K\"ahler 4-Manifolds}
\author{Nikos Georgiou}
\address{Nikos Georgiou\\
          Department of Computing and Mathematics\\
          Institute of Technology, Tralee\\
          Clash\
          Tralee\\
          Co. Kerry\\
          Ireland.}
\email{nikos.georgiou@@research.ittralee.ie}

\author{Brendan Guilfoyle}
\address{Brendan Guilfoyle\\
          Department of Computing and Mathematics \\
          Institute of Technology, Tralee \\
          Clash \\
          Tralee  \\
          Co. Kerry \\
          Ireland.}
\email{brendan.guilfoyle@@ittralee.ie}

\author{Wilhelm Klingenberg}
\address{Wilhelm Klingenberg\\
 Department of Mathematical Sciences\\
 University of Durham\\
 Durham DH1 3LE\\
 United Kingdom.}
\email{wilhelm.klingenberg@@durham.ac.uk }

\date{22nd October 2008}

\begin{abstract}
We study the totally null surfaces of the neutral K\"ahler metric on certain 4-manifolds. The tangent spaces of totally null
surfaces are either self-dual ($\alpha$-planes) or anti-self-dual ($\beta$-planes) and so we consider $\alpha$-surfaces and $\beta$-surfaces.
The metric of the examples we study, which include the spaces of oriented geodesics of 3-manifolds of constant curvature, are anti-self-dual, and so 
it is well-known that the $\alpha$-planes are integrable and $\alpha$-surfaces
exist. These are holomorphic Lagrangian surfaces, which for the geodesic spaces correspond to totally umbilic foliations of the underlying 
3-manifold. 

The $\beta$-surfaces are less known and our interest is mainly in their description. In particular, we classify the $\beta$-surfaces
of the neutral K\"ahler metric on $TN$, the tangent bundle to a Riemannian 2-manifold $N$. These include the spaces of oriented geodesics in Euclidean
and Lorentz 3-space, for which we show that the $\beta$-surfaces are affine tangent bundles to curves of constant geodesic curvature on $S^2$ and $H^2$, 
respectively. In addition, we construct the $\beta$-surfaces of the space of oriented geodesics of hyperbolic 3-space.

\end{abstract}

\maketitle

\section{Introduction}

Neutral K\"ahler 4-manifolds exhibit remarkably different behaviour than their positive-definite counterparts. The 
failure of the complex structure $J$ to tame the symplectic structure $\Omega$ means that 2-planes in the tangent 
space of a point can be both holomorphic and Lagrangian. Under favorable conditions (namely the vanishing of the self-dual
conformal curvature) such planes are integrable and there exist holomorphic Lagrangian surfaces. 

In the space ${\mathbb L}(M)$ of oriented geodesics of a 3-manifold of constant curvature $M$
(on which a natural neutral K\"ahler structure exists) such surfaces play a distinctive role: they correspond to totally 
umbilic foliations of $M$ \cite{gag1} \cite{gak4} \cite{gak5}. 

Holomorphic Lagrangian planes are totally null, that is, the induced metric identically vanishes on the plane. Moreover,
with respect to the Hodge star operator of the neutral metric, the self-dual 2-forms vanish on these planes. There exists 
however another class of totally null planes, upon which the anti-self-dual forms vanish. The former planes are referred to
as $\alpha$-planes, while the latter are $\beta$-planes.

In this note we consider the $\beta$-surfaces in certain neutral K\"ahler 4-manifolds, which include spaces ${\mathbb L}(M)$ of oriented 
geodesics of 3-manifolds $M$ of constant curvature. In the cases of $M={\mathbb E}^3,{\mathbb E}^3_1,{\mathbb H}^3$ we compute the 
$\beta$-surfaces explicitly and show that they include ${\mathbb L}({\mathbb E}^2),{\mathbb L}({\mathbb H}^2)$.

In particular, we prove that

\vspace{0.1in}

\noindent{\bf Main Theorem}:

{\it
A $\beta$-surface in ${\mathbb{L}}({\mathbb{E}}^3)$ is an affine tangent bundle over a curve of constant geodesic curvature in $(S^2,g_{rnd})$.

A $\beta$-surface in ${\mathbb{L}}({\mathbb{E}}_1^3)$ is an affine tangent bundle over a curve of constant geodesic curvature in $(H^2g_{hyp})$.

A $\beta$-surface in ${\mathbb{L}}({\mathbb{H}}^3)$ is a piece of a torus which, up to 
isometry, is either
\begin{enumerate}
\item ${\mathbb{L}}({\mathbb{H}}^2)$, where ${\mathbb{H}}^2\subset{\mathbb{H}}^3$, or
\item ${\cal C}_1\times {{\cal C}}_2\subset\; S^2\times S^2-\bar{\Delta} $, where the ${\cal C}_1$ is a circle given by the
intersection of the 2-sphere and a plane containing the north pole, and ${{\cal C}}_2$ is the image of ${\cal C}_1$ under 
reflection in the horizontal plane through the origin.
\end{enumerate}
}
\vspace{0.1in}

In the next section we discuss self-duality for planes in neutral K\"ahler 4-manifolds and their properties. We then turn 
to the neutral metric on $TN$ and the special case ${\mathbb L}({\mathbb E}^3)$ and ${\mathbb L}({\mathbb E}^3_1)$. In the final section we characterize the
$\beta$-surfaces in ${\mathbb L}({\mathbb H}^3)$.

\vspace{0.2in}

\section{Neutral Metrics on 4-Manifolds}

\subsection{Self-dual and anti-self-dual 2-forms}

Consider the neutral metric $G$ on ${\mathbb R}^4$ given in standard coordinates ($x^1,x^2,x^3,x^4$) by
\[
ds^2=(dx^1)^2+(dx^2)^2-(dx^3)^2-(dx^4)^2.
\]
Throughout, we denote ${\mathbb R}^4$ endowed with this metric by ${\mathbb R}^{2,2}$.

The space of 2-forms on ${\mathbb R}^{2,2}$ is a 6-dimensional linear space that splits naturally with respect to the Hodge
star operator $*$ of $G$ into two 3-dimensional spaces: $\Lambda^2=\Lambda^2_+\oplus\Lambda^2_-$, the space of
self-dual and anti-self-dual 2-forms. Thus, if $\omega\in \Lambda^2$, then $\omega=\omega_++\omega_-$, where
$*\omega_+=\omega_+$ and $*\omega_-=-\omega_-$.

We can easily find a basis for $\Lambda^2_+$ and $\Lambda^2_-$. First, define the {\it double null} basis of 1-forms:
\[
\Theta^1=dx^1+dx^3, \qquad \Theta^2=dx^2- dx^4, \qquad \Theta^3=-dx^2- dx^4, \qquad \Theta^4=dx^1- dx^3, 
\]
so that the metric is
\[
ds^2=\Theta^1\otimes\Theta^4-\Theta^2\otimes\Theta^3.
\]

\vspace{0.1in}

\begin{Prop}\label{p:dual}
If $\omega\in \Lambda^2=\Lambda^2_+\oplus\Lambda^2_-$, with $\omega=\omega_++\omega_-$, then
\[
\omega_+=a_1\Theta^1\wedge\Theta^2+b_1\Theta^3\wedge\Theta^4+c_1(\Theta^1\wedge\Theta^4-\Theta^2\wedge\Theta^3),
\]
\[
\omega_-=a_2\Theta^1\wedge\Theta^3+b_2\Theta^2\wedge\Theta^4+c_2(\Theta^1\wedge\Theta^4+\Theta^2\wedge\Theta^3),
\]
for $a_1,b_1,c_1,a_2,b_2,c_2\in{\mathbb R}$.
\end{Prop}
\begin{pf}
This follows from computing the Hodge star operator acting on 2-forms:
\[
*(\Theta^1\wedge\Theta^4)=-\Theta^2\wedge\Theta^3, \qquad *(\Theta^2\wedge\Theta^4)=-\Theta^2\wedge\Theta^4, \qquad 
*(\Theta^1\wedge\Theta^3)=-\Theta^1\wedge\Theta^3,
\]
\[
*(\Theta^3\wedge\Theta^4)=\Theta^3\wedge\Theta^4, \qquad *(\Theta^1\wedge\Theta^2)=\Theta^1\wedge\Theta^2.
\]
\end{pf}

\vspace{0.2in}

\subsection{Totally Null Planes}

\begin{Def}
A plane $P\subset{\mathbb R}^{2,2}$ is {\it totally null} if every vector in $P$ is null with respect to $G$, and the inner 
product of any two vectors in $P$ is zero. 

A plane  $P$ is {\it self-dual} if $\omega_+(P)=0$ for all 
$\omega_+\in\Lambda^2_+$, and  {\it anti-self-dual} if $\omega_-(P)=0$ for all $\omega_-\in\Lambda^2_-$. 
Self-dual planes are also called $\alpha$-{\it planes}, while anti-self-dual planes are called $\beta$-{\it planes}. 
\end{Def}

\vspace{0.1in}

\begin{Prop}\label{p:tnimpasd}
A plane $P$ is totally null iff $P$ is either self-dual or anti-self-dual.
\end{Prop}
\begin{pf}
Suppose all self-dual forms vanish on $P$ and let $\{V,W\}$ be a basis for $P$. Let ($e_1,e_2,e_3,e_4$) be the vector
basis of ${\mathbb R}^{2,2}$ that is dual to ($\Theta^1,\Theta^2,\Theta^3,\Theta^4$) and $V=V^je_j$, $W=W^je_j$. Since
all of the self-dual 2-forms vanish on $P$, we have from the expression of $\omega_+$ in Proposition \ref{p:dual} that
\begin{equation}\label{e:sdtn1}
V^1W^2=W^1V^2
\qquad\qquad
V^3W^4=W^3V^4,
\end{equation}
\begin{equation}\label{e:sdtn3}
V^1W^4-V^2W^3=W^1V^4-W^2V^3.
\end{equation}
We can assume without loss of generality that $V$ and $W$ are orthogonal: $G(V,W)=0$, which in frame components says that 
\[
V^1W^4+W^1V^4=V^2W^3+W^2V^3.
\]
Combining this with equation (\ref{e:sdtn3}) we have that
\begin{equation}\label{e:sdtn4}
V^1W^4=V^2W^3
\qquad\qquad
W^1V^4=W^2V^3.
\end{equation}
Multiplying the first equation of (\ref{e:sdtn4}) by $W^1$ we have 
\[
V^1W^4W^1=V^2W^3W^1,
\]
which, by virtue of the first equation of (\ref{e:sdtn1}), is
\[
V^1W^4W^1=W^2W^3V^1.
\]
Thus
\[
G(W,W)V^1=2(W^1W^4-W^2W^3)V^1=0.
\]
Similarly, multiplying the first equation of (\ref{e:sdtn4}) by $W^2$, and the second equation by $W^3$ and $W^4$,
applying equations (\ref{e:sdtn1}), we find that
\[
G(W,W)V^2=G(W,W)V^3=G(W,W)V^4=0.
\]
Thus, either $G(W,W)=0$ or $V=0$. Since the latter is not true, we conclude that $W$ is a null vector.

On the other hand,  multiplying the second equation of (\ref{e:sdtn4}) by $V^1$ and $V^2$, and the first by $V^3$ and 
$V^4$, utilizing equations (\ref{e:sdtn1}), we have
\[
G(V,V)W^1=G(V,V)W^2=G(V,V)W^3=G(V,V)W^4=0.
\]
Thus $V$ is also a null vector, and the plane spanned by $V$ and $W$ is totally null, as claimed.

An analogous argument establishes that a plane on which all anti-self-dual 2-forms vanish is totally null.

Conversely, suppose that a plane $P$ is totally null. That is, in terms of a vector basis $V$ and $W$ as before
\begin{equation}\label{e:tnull1}
V^1V^4=V^2V^3 \qquad \qquad W^1W^4=W^2W^3,
\end{equation}
\begin{equation}\label{e:tnull2}
V^1W^4+V^4W^1-V^2W^3-V^3W^2=0.
\end{equation}
Multiplying equation (\ref{e:tnull2}) by $V^1$,$V^3$, $W^1$ and $W^3$, yields, with the aid of equations (\ref{e:tnull1}):
\begin{equation}\label{e:tnull3}
V^2(V^3W^1-V^1W^3)=V^1(V^3W^2-V^1W^4),
\end{equation}
\begin{equation}\label{e:tnull4}
V^4(V^3W^1-V^1W^3)=V^3(V^3W^2-V^1W^4),
\end{equation}
\begin{equation}\label{e:tnull5}
W^2(V^1W^3-V^3W^1)=W^1(V^2W^3-V^4W^1),
\end{equation}
\begin{equation}\label{e:tnull6}
W^4(V^1W^3-V^3W^1)=W^3(V^2W^3-V^4W^1).
\end{equation}
Now, adding $V^1$ times equation (\ref{e:tnull5}), $W^1$ times equation (\ref{e:tnull3}), $V^3$ times equation 
(\ref{e:tnull6}) and $W^3$ times equation (\ref{e:tnull4}) and using equation (\ref{e:tnull2}), we obtain
\begin{equation}\label{e:tnull7}
(V^1W^2-V^2W^1+V^3W^4-V^4W^3)(V^1W^3-V^3W^1)=0.
\end{equation}
By a similar manipulation we find that
\begin{equation}\label{e:tnull8}
(V^1W^2-V^2W^1+V^3W^4-V^4W^3)(V^2W^4-V^4W^2)=0.
\end{equation}

Now suppose that $P$, in addition to being totally null, is Lagrangian. If $J(V)$ is not in $P$, then, since 
$G(W,J(V))=\Omega(W,V)=0$, the metric would be identically zero on the 3-space spanned by $\{V,W,J(V)\}$. For a 
non-degenerate metric $G$ on ${\mathbb R}^{2,2}$ this is not possible. Thus  $J(V)\in P$ and so $P$ is a complex plane.
It follows easily that $P$ is self-dual.

On the other hand, suppose that the totally null plane $P$ is not Lagrangian. Then $\Omega(V,W)\neq0$ or
\[
V^1W^2-V^2W^1+V^3W^4-V^4W^3\neq0.
\]
By equations (\ref{e:tnull7}) and (\ref{e:tnull8}), we have $V^1W^3-V^3W^1=V^2W^4-V^4W^2=0$. Moreover,
substituting these in (\ref{e:tnull3}) to (\ref{e:tnull6}) we conclude that $V^1W^4-V^4W^1+V^2W^3-V^3W^2=0$.
Then, by Proposition \ref{p:dual} we must have $\omega_-(V,W)=0$, which completes the result.

\end{pf}

\vspace{0.1in}

\subsection{K\"ahler Structure on ${\mathbb R}^{2,2}$}

Up to an overall sign, there are two complex structures on ${\mathbb R}^{2,2}$ that are compatible with the metric $G$:
\[
J(X^1,X^2,X^3,X^4)=(-X^2,X^1,-X^4,X^3),
\]
and 
\[
J'(X^1,X^2,X^3,X^4)=(-X^2,X^1,X^4,-X^3).
\]
By compatibility we mean that $G(J\cdot,J\cdot)=G(\cdot,\cdot)$, and similarly for $J'$.

We can utilize these and define two symplectic forms by $\Omega=G(\cdot,J\cdot)$ and $\Omega'=G(\cdot,J'\cdot)$. That
is
\[
\Omega=dx^1\wedge dx^2-dx^3\wedge dx^4 
\qquad\qquad
\Omega'=dx^1\wedge dx^2+dx^3\wedge dx^4.
\]
Thus, the symplectic 2-form $\Omega$ is self-dual while $\Omega'$ is anti-self-dual.

Moreover

\begin{Prop}
An $\alpha$-plane is holomorphic and Lagrangian with respect to $(J,\Omega)$, while a $\beta$-plane is
holomorphic and Lagrangian with respect to $(J',\Omega')$.
\end{Prop}
\begin{pf}
The proof follows from arguments similar to those of Proposition \ref{p:tnimpasd}.
\end{pf}

Given a null vector $V$ in ${\mathbb R}^{2,2}$, the planes spanned by $\{V,J(V)\}$ and $\{V,J'(V)\}$ are 
easily seen to be totally null. More explicitly, the set of totally null planes is, in fact, $S^1\cup S^1$, 
which can be parameterized as follows.
For $a,b\in{\mathbb R}$, $\phi\in[0,2\pi)$ and $\epsilon=\pm 1$, consider the vector in ${\mathbb R}^{2,2}$
given by
\[
V^\epsilon_{\phi}(a,b)=\left(a\cos\phi+b\sin\phi,a\sin\phi-b\cos\phi,
                          a,-\epsilon b\right).
\]
Let $P^\epsilon_{\phi}$ be the plane containing $V^\epsilon_{\phi}(a,b)$ as $a$ and $b$
vary over ${\mathbb R}$. Then a quick check shows that $P^+_\phi$ is self-dual, while $P^-_\phi$ is anti-self-dual.

\subsection{Neutral K\"ahler 4-manifolds}

Let $(M,G,J,\Omega)$ be a smooth neutral K\"ahler 4-manifold. Thus $M$ is a smooth 4-manifold, $G$ is a neutral metric, while $J$
is a complex structure that is compatible with $G$ and $\Omega(\cdot,\cdot)=G(J\cdot,\cdot)$ is a closed non-degenerate (symplectic) 2-form. 

The existence of a unitary frame at a point of $M$ implies that it is possible to apply the algebra of the last section pointwise on $M$, and 
we therefore have $S^1\cup S^1$ worth of totally null planes at each point. On a compact 4-manifold, the existence of an oriented 2-dimensional 
distribution implies topological restrictions on $M$ \cite{HaH}, and so not every compact 4-manifold admits a neutral K\"ahler structure. 
However, the examples we consider are non-compact and the neutral K\"ahler structure will be given explicitly. 

On any (pseudo)-Riemannian 4-manifold $(M,G)$ the Riemann curvature tensor can be considered as an endomorphism of 
$\Lambda^2(M)$. The splitting $\Lambda^2(M)=\Lambda^2_+(M)\oplus\Lambda^2_-(M)$ with respect to the Hodge star operator $*$ 
yields a block decomposition of the Riemann curvature tensor
\[
{\mbox {Riem}}=\left[\begin{matrix}
&   & & & \\
& {\mbox {Weyl}}^++{\textstyle{\frac{1}{12}}}{\mbox {R}} & & {\mbox {Ric}} & \\ 
&   & & & \\
&  {\mbox {Ric}} & & {\mbox {Weyl}}^-+{\textstyle{\frac{1}{12}}}{\mbox {R}} &\\
&   & & & \\
\end{matrix}\right],
\]
where Ric is the Ricci tensor, R is the scalar curvature and Weyl$^{\pm}$ are the self- and anti-self-dual Weyl curvature
tensors \cite{besse}.

\begin{Def}
A (pseudo)-Riemannian 4-manifold ($M,G$) is {\it anti-self-dual} if the self-dual part of the Weyl conformal curvature 
tensor vanishes: Weyl$^+=0$.
\end{Def}

\vspace{0.1in}

A well-known result of Penrose states:

\begin{Thm} \cite{daw}\cite{penrose}\label{t:pen}
The $\alpha$-surfaces of a neutral K\"ahler 4-manifold ($M,G$) are integrable iff ($M,G$) is anti-self-dual.
\end{Thm}

\vspace{0.2in}

\section{Neutral K\"ahler Metric on $TN$}

Let $(N,g)$ be a Riemannian 2-manifold and consider the total space $TN$ of the tangent bundle to $N$.
Choose conformal coordinates $\xi$ on $N$ so that $ds^2=e^{2u}d\xi d\bar{\xi}$ for some function $u=u(\xi,\bar{\xi})$, 
and the corresponding complex coordinates $(\xi,\eta)$ on $TN$ obtained by identifying
\[
(\xi,\eta) \leftrightarrow \eta\frac{\partial}{\partial \xi}+\bar{\eta}\frac{\partial}{\partial \bar{\xi}}
                   \in T_\xi N.
\] 
The coordinates $(\xi,\eta)$ define a natural complex structure on $TN$ by
\[
{\mathbb J}\left(\frac{\partial}{\partial\xi}\right)=i\frac{\partial}{\partial\xi}
\qquad\qquad
{\mathbb J}\left(\frac{\partial}{\partial\eta}\right)=i\frac{\partial}{\partial\eta}.
\]
In \cite{gak4} a neutral K\"ahler structure was introduced on $TN$. In the above coordinate system, the symplectic 2-form is
\begin{equation}\label{e:symp}
\Omega=2e^{2u}{\Bbb{R}}\mbox{e}\left(d\eta\wedge d\bar{\xi}+2\eta\partial_\xi u\; d\xi\wedge d\bar{\xi}\right),
\end{equation}
while the neutral metric ${\Bbb{G}}$ is 
\begin{equation}\label{e:metric1}
{\Bbb{G}}=2e^{2u}{\Bbb{I}}\mbox{m}\left(d\bar{\eta} d\xi-2\eta\partial_\xi u \;d\xi d\bar{\xi}\right).
\end{equation}
Here we have introduced the notation $\partial_\xi$ for differentiation with respect to $\xi$.

\begin{Note}
When $u=0$ we retrieve the neutral K\"ahler metric on ${\mathbb R}^4={\mbox{T}}{\mathbb R}^2$, where
\[
\xi={\textstyle{\frac{1}{2}}}\left[x^1+x^3+i(x^2+x^4)\right]
\qquad\qquad
\eta={\textstyle{\frac{1}{2}}}\left[x^2-x^4+i(-x^1+x^3)\right],
\]
or
\[
x^1={\textstyle{\frac{1}{2}}}\left[\xi+\bar{\xi}+i(\eta-\bar{\eta})\right]
\qquad\qquad
x^2={\textstyle{\frac{1}{2}}}\left[-i(\xi-\bar{\xi})+\eta+\bar{\eta}\right],
\]
\[
x^3={\textstyle{\frac{1}{2}}}\left[\xi+\bar{\xi}-i(\eta-\bar{\eta})\right]
\qquad\qquad
x^4={\textstyle{\frac{1}{2}}}\left[-i(\xi-\bar{\xi})-\eta-\bar{\eta}\right].
\]
\end{Note}

\vspace{0.1in}

\begin{Prop}
The double null basis for $(TN,{\mathbb G}$) is
\[
\Theta^1=2{\Bbb{R}}\mbox{e} (d\xi)
\qquad\qquad
\Theta^2=2e^{2u}{\Bbb{R}}\mbox{e} \left(d\eta+2\eta\partial_\xi u\; d\xi\right),
\]
\[
\Theta^3=2{\Bbb{I}}\mbox{m} (d\xi)
\qquad\qquad
\Theta^4=2e^{2u}{\Bbb{I}}\mbox{m} \left(d\eta+2\eta\partial_\xi u\; d\xi\right).
\]
\end{Prop}
\begin{pf}
A straight-forward check shows that
\[
ds^2=\Theta^1\otimes\Theta^4-\Theta^2\otimes\Theta^3,
\]
as claimed.
\end{pf}

\vspace{0.1in}

The coordinate expressions for self-dual and anti-self-dual 2-forms on $TN$ are

\begin{Prop}\label{p:dualTN}
If $\omega\in \Lambda^2(TN)=\Lambda^2_+(TN)\oplus\Lambda^2_-(TN)$, with $\omega=\omega_++\omega_-$, then
\begin{align}
\omega_+&=a_1(d\xi\wedge d\eta+d\bar{\xi}\wedge d\bar{\eta})+b_1[d\xi\wedge d\bar{\eta}+d\bar{\xi}\wedge d\eta
     +2(\bar{\eta}\partial_{\bar{\xi}}u-\eta\partial_\xi u)d\xi\wedge d\bar{\xi}]\nonumber \\
      &\qquad\qquad\qquad+ic_1(d\xi\wedge d\eta-d\bar{\xi}\wedge d\bar{\eta}),\nonumber
\end{align}
\begin{align}
\omega_-&=ia_2d\xi\wedge d\bar{\xi}+ib_2[d\xi\wedge d\bar{\eta}-d\bar{\xi}\wedge d\eta
     +2(\bar{\eta}\partial_{\bar{\xi}}u+\eta\partial_\xi u)d\xi\wedge d\bar{\xi}]\nonumber \\
      &\qquad\qquad+ic_2(d\eta\wedge d\bar{\eta}+2\eta\partial_\xi u d\xi\wedge d\bar{\eta}
          +2\bar{\eta}\partial_{\bar{\xi}} u d\bar{\xi}\wedge d\eta
           +4\eta\bar{\eta}\partial_\xi u\partial_{\bar{\xi}} ud\xi\wedge d\bar{\xi}),\nonumber
\end{align}
for $a_1,b_1,c_1,a_2,b_2,c_2\in{\mathbb R}$.
\end{Prop}

\vspace{0.1in}

\subsection{$\alpha$-surfaces in $TN$}

We first note that 

\begin{Prop}
The neutral K\"ahler metric ${\mathbb G}$ on $TN$ is anti-self-dual.
\end{Prop}
\begin{pf}
A calculation using the coordinate expression (\ref{e:metric1}) of the metric shows that the only non-vanishing 
component of the conformal curvature tensor is
\[
W_{\xi\bar{\xi}}^{\;\;\;\;\eta\bar{\eta}}=i(\eta\partial_\xi \kappa-\bar{\eta}\partial_{\bar{\xi}}\kappa),
\]
where $\kappa$ is the Gauss curvature of $(N,g)$.

Thus, from Proposition \ref{p:dualTN}, for any $\omega_+\in\Lambda^2_+(TN)$, $W(\omega_+)=0$. That is, the metric
is anti-self-dual.
\end{pf}
\vspace{0.1in}

By applying Theorem \ref{t:pen} we have

\begin{Cor}
There exists $\alpha$-surfaces, i.e. holomorphic Lagrangian surfaces, in $(TN,{\mathbb J},\Omega)$. 
\end{Cor}

\vspace{0.1in}

\subsection{$\beta$-surfaces in $TN$}

\begin{Prop}\label{p:betaTN}
An immersed surface  $\Sigma\subset TN$ is a $\beta$-{\it surface} iff locally it is given by $(s,t)\rightarrow(\xi(s,t),\eta(s,t))$
where
\[
\xi=se^{iC_0}+\xi_0 \qquad\qquad
\eta=(te^{iC_0}+\eta_0)e^{-2u},
\]
for $C_0\in{\mathbb R}$ and $\xi_0,\eta_0\in{\mathbb C}$.
\end{Prop}
\begin{pf}
By Proposition \ref{p:dualTN} surface $f:\Sigma\rightarrow TN$ is a $\beta$-{\it surface} iff
\begin{equation}\label{e:beta11}
f^*(d\xi\wedge d\bar{\xi})=0
\qquad\qquad
f^*(d(\eta e^{2u})\wedge d(\bar{\eta} e^{2u}))=0,
\end{equation}
and
\begin{equation}\label{e:beta22}
f^*(d\xi\wedge d(\bar{\eta} e^{2u})-d\bar{\xi}\wedge d(\eta e^{2u}))=0.
\end{equation}

The first equation of (\ref{e:beta11}) implies that the map $(s,t)\rightarrow \xi(s,t)$ is not of maximal rank, and as it
cannot be of rank zero (as this would mean that $\Sigma$ is a fibre of $\pi:TN\rightarrow N$, and is therefore an 
$\alpha$-surface) it must be of rank 1. By the implicit function theorem either
\[
\xi(s,t)=\xi(s,t(s)) \qquad {\mbox {or}} \qquad \xi(s,t)=\xi(s(t),t).
\]
Without loss of generality, we will asume the former: $\xi=\xi(s)$.

Similarly, the second equation of (\ref{e:beta11}) implies that either
\[
\eta e^{2u}=\psi(s,t)=\psi(s,t(s)) \qquad {\mbox {or}} \qquad \eta e^{2u}=\psi(s,t)=\psi(s(t),t).
\]
Here, we must have the latter $\eta e^{2u}=\psi(t)$, or else the surface $\Sigma$ would be singular.

Turning now to equation of (\ref{e:beta22}), we have
\[
\frac{d \xi}{d s}\frac{d \bar{\psi}}{d t}=\frac{d \bar{\xi}}{d s}\frac{d \psi}{d t}.
\]
By separation of variables we see that
\[
\frac{d \xi}{d s}=e^{2iC_0}\frac{d \bar{\xi}}{d s}
\qquad\qquad
\frac{d \psi}{d s}=e^{2iC_0}\frac{d \bar{\psi}}{d s},
\]
for some real constant $C_0$. These can be integrated to
\[
\xi=h_1(s)e^{iC_0}+\xi_0 \qquad \qquad \eta=(h_2(t)e^{iC_0}+\eta_0)e^{-2u},
\]
for complex constants $\xi_0$ and $\eta_0$ and real functions $h_1$ and $h_2$ of $s$ and $t$, respectively.

Finally, we can reparameterize $s$ and $t$ so that $h_1=s$ and $h_2=t$, as claimed.
\end{pf}

\vspace{0.1in}

\subsection{The Oriented Geodesic Spaces $TS^2$ and $TH^2$}

In the cases where $N=S^2$ or $N=H^2$ endowed with a metric of constant Gauss curvature 
($e^{2u}=4(1\pm\xi\bar{\xi})^{-2}$), the above construction yields the
neutral K\"ahler metric on the space ${\mathbb{L}}({\mathbb{E}}^{3})$ of oriented affine lines or on the space 
${\mathbb{L}}({\mathbb{E}}^{3}_{1})$ of 
future-pointing time-like lines, in ${\mathbb{E}}^3$ or ${\mathbb{E}}^3_1$ (respectively) \cite{gak5}. 

In what follows we consider only the Euclidean case, although analogous results hold for the Lorentz case.
We define the map $\Phi$ which sends ${\mathbb{L}}({\mathbb{E}}^{3})\times{\mathbb{R}}$ 
to ${\mathbb{E}}^3$ as follows: $\Phi$ takes an oriented line $\gamma$
and a real number $r$ to that point in ${\mathbb{E}}^3$ which lies on $\gamma$ and is an affine 
parameter distance $r$ from the point on $\gamma$ closest to the origin.

\vspace{0.1in}

\begin{Prop}\cite{gak4}
The map can be written as $\Phi((\xi,\eta),r)=(z,t)\in{\mathbb{C}}\oplus{\mathbb{R}}={\mathbb{E}}^3$  
where the local coordinate expressions are:
\begin{equation}\label{e:coord1}
z=\frac{2(\eta-\bar{\eta}\xi^2)+2\xi(1+\xi\bar{\xi})r}{(1+\xi\bar{\xi})^2}
\qquad\qquad
t=\frac{- 2(\eta\bar{\xi}+\bar{\eta}\xi)+(1-\xi^2\bar{\xi}^2)r}{(1+\xi\bar{\xi})^2},
\end{equation}
and
\begin{equation}\label{e:coord2}
\eta={\textstyle{\frac{1}{2}}}(z-2t\xi-\bar{z}\xi^2) \qquad\qquad 
          r=\frac{\bar{\xi}z+\xi\bar{z}+(1-\xi\bar{\xi})t}{1+\xi\bar{\xi}}.
\end{equation}

\end{Prop}

\vspace{0.1in}

For $\alpha$-surfaces, we have

\begin{Prop}
A holomorphic Lagrangian surface in $TS^2$ corresponds to the oriented normals to totally umbilic surfaces in ${\mathbb E}^3$ i.e. round spheres or
planes.
\end{Prop}

On the other hand:

\begin{Prop}
A $\beta$-surface in $TS^2$ is an affine tangent bundle over a curve of constant geodesic curvature in $(S^2,g_{rnd})$.
\end{Prop}

\begin{pf}
By Proposition \ref{p:betaTN}, the $\beta$-surfaces are given by
\[
\xi=se^{iC_0}+\xi_0 \qquad\qquad
\eta=(1+\xi\bar{\xi})^2(te^{iC_0}+\eta_0).
\]
Clearly this is a real line bundle over a curve on $S^2$. 
By a rotation this can be simplified to 
\[
\xi=s+\xi_0e^{-iC_0}\qquad\qquad
\eta=(1+\xi\bar{\xi})^2(t+\eta_0e^{-iC_0}),
\]
and after an affine reparameterization of $s$ and $t$ we can set
\[
\xi=s+iC_1\qquad\qquad
\eta=(1+\xi\bar{\xi})^2(t+iC_2).
\]
Projecting onto $S^2$ we get the curve $\xi=s+iC_1$ with unit tangent $\vec{T}$ and normal vector $\vec{N}$ (with respect to the round
metric)
\[
\vec{T}=\frac{(1+\xi\bar{\xi})}{2\sqrt{2}}\left(\frac{\partial}{\partial \xi}+\frac{\partial}{\partial \bar{\xi}}\right)
\qquad\qquad
\vec{N}=\frac{i(1+\xi\bar{\xi})}{2\sqrt{2}}\left(\frac{\partial}{\partial \xi}-\frac{\partial}{\partial \bar{\xi}}\right).
\]
Considered as a set of vectors on $S^2$, the $\beta$-surface is
\begin{align}
\eta\frac{\partial}{\partial \xi}+\bar{\eta}\frac{\partial}{\partial \bar{\xi}}
   &=(1+\xi\bar{\xi})^2(t+iC_2)\frac{\partial}{\partial \xi}+(1+\xi\bar{\xi})^2(t-iC_2)\frac{\partial}{\partial \bar{\xi}}\nonumber\\
   &=2\sqrt{2}(1+\xi\bar{\xi})(t\vec{T}+C_2\vec{N}).\nonumber
\end{align}
These form a real line bundle over the base curve - which do not pass through the origin in the fibre of $TS^2$ for $C_2\neq0$. For $C_2=0$,
this is exactly the tangent bundle to the curve.

The geodesic curvature of this curve is
\begin{align}
g(\vec{N},\nabla_{\vec{T}} \vec{T})&=N_kT^j(\partial_jT^k+\Gamma_{jl}^kT^l)\nonumber\\
&=N_kT^j\partial_jT^k+N^kT^jT^l(2\partial_jg_{lk}-\partial_kg_{jl})\nonumber\\
&=\sqrt{2}C_1.\nonumber
\end{align}
\end{pf}

\vspace{0.1in}

A similar calculation establishes:

\begin{Prop}
A $\beta$-surface in $TH^2$ is an affine tangent bundle over a curve of constant geodesic curvature in $(H^2,g_{hyp})$.
\end{Prop}

\vspace{0.1in}

We also have the following:

\begin{Cor}
Given an affine plane $P$ in ${\mathbb E}^3$, the set ${\mathbb L}({\mathbb E}^2)$ of oriented lines contained in 
$P$ is a  $\beta$-surface in $TS^2$.
\end{Cor}
\begin{pf}

By Proposition \ref{p:betaTN}, the $\beta$-surfaces are given by
\[
\xi=se^{iC_0}+\xi_0 \qquad\qquad
\eta=(1+\xi\bar{\xi})^2(te^{iC_0}+\eta_0)
\]
Isometries of ${\mathbb E}^3$ induce isometries on $TS^2$ and hence preserves 
$\beta$-surfaces. Thus we can translate and rotate $P$ so that it is vertical and contains the $t$-axis.
Thus we can consider the $\beta$-surface $\Sigma$ with $\xi_0=\eta_0=0$, and then using the map $\Phi$ we find the two parameter 
family of oriented lines in ${\mathbb{E}}^3$ to be
\[
z=\frac{2[(1-s^4)t+sr]}{1+ s^2}e^{iC_0} \qquad
t=\frac{- 4s(1+ s^2)t+(1-s^2)r}{1+ s^2}
\]
This is a vertical plane containing the $t$-axis, and $\Sigma$ consists of all the oriented lines in this plane.
\end{pf}

\vspace{0.2in}

\section{Oriented Geodesics in Hyperbolic 3-space}

We briefly recall the basic construction of the canonical neutral K\"ahler metric on the space 
${\mathbb{L}}({\mathbb{H}}^3)$ of oriented geodesics of ${\mathbb{H}}^3$ - further details can be found in \cite{gag1}.

Consider the 4-manifold ${\mathbb P}^1\times{\mathbb P}^1$ endowed with the canonical complex structure
${\mathbb J}=j\oplus j$ and complex coordinates $\mu_1$ and $\mu_2$. If we let 
$\overline{\Delta}=\{(\mu_1,\mu_2):\mu_1\bar{\mu}_2=-1\}$ then 
${\mathbb L}({\mathbb H}^3)={\mathbb P}^1\times{\mathbb P}^1-\overline{\Delta}$. We introduce the neutral K\"ahler metric and symplectic
form on ${\mathbb L}({\mathbb H}^3)$ by
\begin{equation}\label{e:metric2}
{\mathbb{G}}=-i\left[\frac{1}{(1+\mu_{1}\bar{\mu}_{2})^2}d\mu_{1}\otimes d\bar{\mu}_{2}
                   -\frac{1}{(1+\bar{\mu}_{1}\mu_{2})^2}d\bar{\mu}_{1}\otimes d\mu_{2}\right],
\end{equation}
and
\begin{equation}\label{e:sympl}
\Omega=-\left[\frac{1}{(1+\mu_{1}\bar{\mu}_{2})^2}d\mu_{1}\wedge d\bar{\mu}_{2}
           +\frac{1}{(1+\bar{\mu}_{1}\mu_{2})^2}d\bar{\mu}_{1}\wedge d\mu_{2}\right].
\end{equation}

\vspace{0.1in}

\begin{Prop}
A double null basis for $({\mathbb L}({\mathbb H}^3),\; {\mathbb{G}})$ is 
\[
\Theta^1={\mathbb{R}}e\left(\frac{d\mu_1}{1+\mu_1\bar{\mu}_2}-\frac{d\mu_2}{1+\bar{\mu}_1\mu_2}\right)\qquad\qquad \Theta^2={\mathbb{R}}e\left(\frac{d\mu_1}{1+\mu_1\bar{\mu}_2}+\frac{d\mu_2}{1+\bar{\mu}_1\mu_2}\right),
\]
\[
\Theta^3=-{\mathbb{I}}m\left(\frac{d\mu_1}{1+\mu_1\bar{\mu}_2}-\frac{d\mu_2}{1+\bar{\mu}_1\mu_2}\right)\qquad\qquad \Theta^4=-{\mathbb{I}}m\left(\frac{d\mu_1}{1+\mu_1\bar{\mu}_2}+\frac{d\mu_2}{1+\bar{\mu}_1\mu_2}\right).
\]
\end{Prop}
\begin{pf}
A straight-forward computation shows that
\[
ds^2=\Theta^1\otimes\Theta^4-\Theta^2\otimes\Theta^3,
\]
as claimed
\end{pf} 
The coordinate expressions for self-dual and anti-self-dual 2 forms on ${\mathbb L}({\mathbb H}^3)$ are easily found to be:
\begin{Prop}
If $\omega\in \Lambda^2 ({\mathbb L}({\mathbb H}^3))=\Lambda_{+}^2({\mathbb L}({\mathbb H}^3))\oplus\Lambda_{-}^2({\mathbb L}({\mathbb H}^3))$, with $\omega=\omega_{+}+\omega_{-}$, then
\[
\omega_{+}=(a_1+ic_1)\frac{d\mu_1\wedge d\mu_2}{|1+\bar{\mu}_1\mu_2|^2}+(a_1-ic_1)\frac{d\bar{\mu}_1\wedge d\bar{\mu}_2}{|1+\bar{\mu}_1\mu_2|^2}+b_1\left[\frac{d\mu_1\wedge d\bar{\mu}_2}{(1+\mu_1\bar{\mu}_2)^2}+\frac{d\bar{\mu}_1\wedge d\mu_2}{(1+\bar{\mu}_1\mu_2)^2}\right],                   
\]
\[
\omega_{-}=-i(a_2+c_2)\frac{d\mu_1\wedge d\bar{\mu}_1}{|1+\bar{\mu}_1\mu_2|^2}-i(a_2-c_2)\frac{d\mu_2\wedge d\bar{\mu}_2}{|1+\bar{\mu}_1\mu_2|^2}+ib_2\left[\frac{d\mu_1\wedge d\bar{\mu}_2}{(1+\mu_1\bar{\mu}_2)^2}-\frac{d\bar{\mu}_1\wedge d\mu_2}{(1+\bar{\mu}_1\mu_2)^2}\right].
\]
\end{Prop}

\vspace{0.1in}

\subsection{$\alpha$-surfaces in ${\mathbb L}({\mathbb H}^3)$}

Once again, the neutral metric on ${\mathbb L}({\mathbb H}^3)$ is anti-self-dual, indeed, it is conformally flat, and so there exists
$\alpha$-surfaces in ${\mathbb L}({\mathbb H}^3)$. These are found to be the normal congruence to the totally umbilic surfaces in ${\mathbb H}^3$:

\begin{Prop}\cite{gag2}
A smooth surface $\Sigma$ in ${\mathbb L}({\mathbb H}^3)$ is totally null iff $\Sigma$ is the oriented normal congruence 
of 
\begin{enumerate}
\item a geodesic sphere, or 
\item a horosphere, or
\item a totally geodesic surface
\end{enumerate}
in ${\mathbb{H}}^3$.
\end{Prop}

\vspace{0.1in}

\subsection{$\beta$-surfaces in ${\mathbb L}({\mathbb H}^3)$}

\begin{Prop}
Let $\Sigma$ be a $\beta$-surface in ${\mathbb{L}}({\mathbb{H}}^3)$. Then $\Sigma$ is a piece of a torus which, up to 
isometry, is either
\begin{enumerate}
\item ${\mathbb{L}}({\mathbb{H}}^2)$, where ${\mathbb{H}}^2\subset{\mathbb{H}}^3$, or
\item ${\cal C}_1\times {{\cal C}}_2\subset\; S^2\times S^2-\bar{\Delta} $, where the ${\cal C}_1$ is a circle given by the
intersection of the 2-sphere and a plane containing the north pole, and ${{\cal C}}_2$ is the image of ${\cal C}_1$ under 
reflection in the horizontal plane through the origin.
\end{enumerate}
\end{Prop}
\begin{pf}

Let $f:\Sigma\rightarrow {\mathbb L}({\mathbb H}^3)$ be an immersed $\beta$-surface. Then for every anti-self-dual 2-form $\omega_{-}$ we have $f^{*}\omega_{-}=0$. Then we obtain the following equations
\begin{equation}\label{e:exisosi1}
f^{*}(d\mu_1\wedge d\bar{\mu}_1)=0
\qquad\qquad
f^{*}(d\mu_2\wedge d\bar{\mu}_2)=0,
\end{equation}
\begin{equation}\label{e:exisosi3}
f^{*}\left(\frac{d\mu_1\wedge d\bar{\mu}_2}{(1+\mu_1\bar{\mu}_2)^2}-\frac{d\bar{\mu}_1\wedge d\mu_2}{(1+\bar{\mu}_1\mu_2)^2}\right)=0.
\end{equation}
The first equation of (\ref{e:exisosi1}) implies that the map $(u,v)\mapsto \mu_1(u,v)$ is not of maximal rank and since it cannot be of rank zero (otherwise $\Sigma$ would be an $\alpha$-surface) it must be of rank 1. By the implicit function theorem either
\[
\mu_1(u,v)=\mu_1(u,v(u))\qquad {\it or}\qquad \mu_1(u,v)=\mu_1(u(v),v).
\]
Without loss of generality, we will assume the former: $\mu_1=\mu_1(u)$.

Similarly, the second equation of (\ref{e:exisosi1}) implies that 
\[
\mu_2(u,v)=\mu_2(u,v(u))\qquad {\it or}\qquad \mu_2(u,v)=\mu_2(u(v),v).
\]
Here, we must have $\mu_2=\mu_2(v)$, or else the surface $\Sigma$ would be singular.

The equation (\ref{e:exisosi3}) yields
\begin{equation}\label{e:exisosi4}
\ln\mu_2-\ln\bar{\mu}_2+\ln(1+\bar{\mu}_1\mu_2)-\ln(1+\mu_1\bar{\mu}_2)=h_1(u)+h_2(v),
\end{equation}
\begin{equation}\label{e:exisosi5}
\ln\bar{\mu}_1-\ln\mu_1+\ln(1+\bar{\mu}_1\mu_2)-\ln(1+\mu_1\bar{\mu}_2)=w_1(u)+w_2(v),
\end{equation}
for some complex functions $h_1,h_2,w_1,w_2$.

If $h_{i}=a_{i}e^{i\phi_{i}}$ for $i=1,2$, where $a_1=a_1(u),\; \phi_1=\phi_1(u)$ and $a_2=a_2(v),\; \phi_2=\phi_2(v)$ 
are real functions, we obtain 
\[
h_1(u)=ia_1\qquad\qquad h_2(v)=ia_2.
\]
By a similar argument, there are real functions $b_1=b_1(u)$ and $b_2=b_2(v)$ such that (\ref{e:exisosi4}) and (\ref{e:exisosi5}) become
\begin{equation}\label{e:exisosi6}
\ln\mu_2-\ln\bar{\mu}_2+\ln(1+\bar{\mu}_1\mu_2)-\ln(1+\mu_1\bar{\mu}_2)=i(a_1(u)+a_2(v)),
\end{equation}
\begin{equation}\label{e:exisosi7}
\ln\bar{\mu}_1-\ln\mu_1+\ln(1+\bar{\mu}_1\mu_2)-\ln(1+\mu_1\bar{\mu}_2)=i(b_1(u)+b_2(v)).
\end{equation}

Finally from combining equations (\ref{e:exisosi6}) and (\ref{e:exisosi7}) we have
\[
\ln\left(\frac{1+\bar{\mu}_1\mu_2}{1+\mu_1\bar{\mu}_2}\right)=-2i(f(u)+g(v)).
\]

We are thus led to consider the curves ${\cal C}_1,{\cal C}_2$ on $S^2$ given locally by non-constant functions 
$\mu_1:{\mathbb R}\rightarrow S^2: u\mapsto \mu_1(u)$ 
and $\mu_2:{\mathbb R}\rightarrow S^2: v\mapsto \mu_2(v)$ which satisfy
\[
1+\mu_1\bar{\mu}_2=(1+\bar{\mu}_1\mu_2)e^{2i(f+g)},
\]
for $f=f(u)$ and $g=g(v)$.

If we switch to polar coordinates  $\mu_1=\lambda_1(u)e^{i\theta_1(u)}$ and $\mu_2=\lambda_2(v)e^{i\theta_2(v)}$, 
this reduces to
\begin{equation}\label{e:beta1}
\sin[f(u)+g(v)]=\lambda_1(u)\lambda_2(v)\sin[\theta_1(u)-f(u)-\theta_2(v)-g(v)].
\end{equation}

By a rotation we can set $\mu_2$ to zero for some $v=v_0$, that is, $\lambda_2(v_0)=0$.  We 
find from equation (\ref{e:beta1}) that
\[
\sin[f(u)+g(v_0)]=0,
\] 
and so letting $g_0=g(v_0)$, we conclude that $f=-g_0$. Putting this back into (\ref{e:beta1}) we have
\begin{equation}\label{e:beta1a}
\sin[g(v)-g_0]=\lambda_1(u)\lambda_2(v)\sin[\theta_1(u)-\theta_2(v)-g(v)+g_0].
\end{equation}
Thus for a fixed $u=u_0$ we have 
\[
\lambda_1(u_0)\lambda_2(v)\sin[\theta_1(u_0)-\theta_2(v)-g(v)+g_0]=\lambda_1(u)\lambda_2(v)\sin[\theta_1(u)-\theta_2(v)-g(v)+g_0],
\]
or, for $v\neq v_0$
\begin{equation}\label{e:beta2}
\lambda_1(u_0)\sin[\theta_1(u_0)-\theta_2(v)-g(v)+g_0]=\lambda_1(u)\sin[\theta_1(u)-\theta_2(v)-g(v)+g_0].
\end{equation}
Differentiating this relationship with respect to $v$ yields
\begin{align}
&\lambda_1(u_0)\cos[\theta_1(u_0)-\theta_2(v)-g(v)+g_0]\;\partial_v(\theta_2+g) \nonumber \\
&\qquad\qquad =\lambda_1(u)\cos[\theta_1(u)-\theta_2(v)-g(v)+g_0]\;\partial_v(\theta_2+g) \label{e:beta3}.
\end{align}

If $\partial_v(\theta_2+g)\neq 0$, then we can cancel this factor and square both sides of equations (\ref{e:beta2}) and 
(\ref{e:beta3}) to find that $\lambda_1=\lambda_1(u_0)$. However, from the functional relation  in equation 
(\ref{e:beta1a}), this means that
$\theta_1$ is also constant. Thus $\mu_1$ would be constant, which is not true.

We conclude that $\partial_v(\theta_2+g)=0$, or equivalently, $g(v)=-\theta_2(v)+g_1$. Substituting this back into
equation (\ref{e:beta1a}) we have
\[
\sin[\theta_2(v)+C_0]=\lambda_1(u)\lambda_2(v)\sin[\theta_1(u)+C_0],
\]
where $C_0=g_0-g_1$.

One solution of this equation is $\theta_1=\theta_2=-C_0$, which is the case 
$\Sigma={\mathbb{L}}({\mathbb{H}}^2)$, where ${\mathbb{H}}^2\subset{\mathbb{H}}^3$. 
Otherwise, we can separate variables
\[
\frac{\sin[\theta_2(v)+C_0)]}{\lambda_2(v)}=\lambda_1(u)\sin[\theta_1(u)+C_0]=C_1\neq0.
\]
This yields
\[
\mu_1=\frac{C_1 e^{i\theta_1(u)}}{\sin[\theta_1(u)+C_0]}
\qquad\qquad
\mu_2=\frac{\sin[\theta_2(v)+C_0]e^{i\theta_2(v)}}{C_1}.
\]
By a rotation of $S^2$ we can set $C_0$ to zero, and with a natural choice of parameterization of the curves, the final 
form is
\[
\mu_1=\frac{C_1 e^{iu}}{\sin u}
\qquad\qquad
\mu_2=\frac{\sin v \;e^{iv}}{C_1},
\]
for $u,v\in[0,2\pi)$.

These are the tori of part (2) in the statement. To see that they are circles note that if we view $S^2$ in
${\mathbb R}^3$ given by
\[
x=\frac{\mu+\bar{\mu}}{1+\mu\bar{\mu}}
\qquad
y=\frac{-i(\mu-\bar{\mu})}{1+\mu\bar{\mu}}
\qquad
z=\frac{1-\mu\bar{\mu}}{1+\mu\bar{\mu}},
\]
then the first curve parameterizes the intersection of $S^2$ with the plane $y+C_1(z-1)=0$, while the second is the 
intersection with the plane $y-C_1(z+1)=0$.

\end{pf}

In the ball model of ${\mathbb H}^3$ these 2-parameter families of geodesics can be visualized as the set of  
geodesics that begin on a circle on the boundary and end on another circle of the same radius on the boundary, the two circles having
a single point of intersection, as illustrated below.

\vspace{0.1in}
\setlength{\epsfxsize}{4.5in}
\begin{center}
   \mbox{\epsfbox{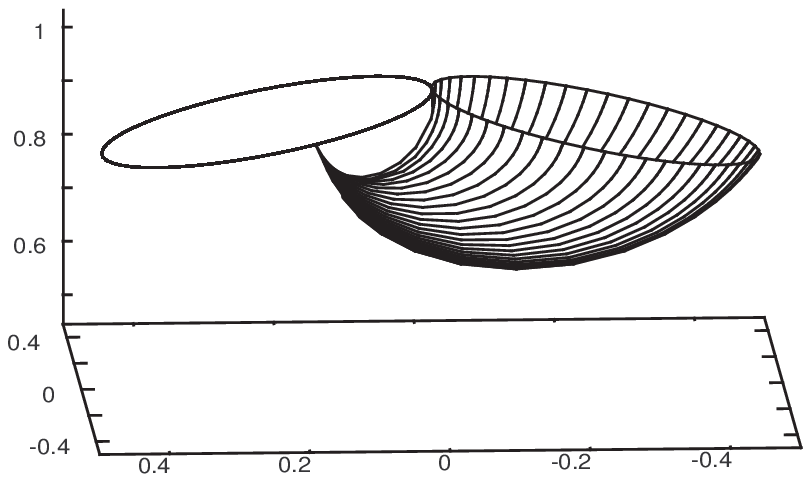}}
\end{center}

\end{document}